\documentclass[12pt]{amsart}
\usepackage{amsmath}
\usepackage{amsfonts}
\usepackage{amsthm}
\usepackage{amssymb}

\theoremstyle{plain}

\newtheorem{Theorem}{Theorem}

\newtheorem{Proposition}[Theorem]{Proposition}

\newtheorem{Question}[Theorem]{Question}

\theoremstyle{remark}

\begin{document}

\title[Combinatorial identities]
{Combinatorial identities related to  eigenfunction decompositions of
Hill operators: Open Questions}

{\author{Plamen Djakov}}


\author{Boris Mityagin}

\address{Sabanci University, Orhanli,
34956 Tuzla, Istanbul, Turkey}
 \email{djakov@sabanciuniv.edu}

\address{Department of Mathematics,
The Ohio State University,
 231 West 18th Ave,
Columbus, OH 43210, USA} \email{mityagin.1@osu.edu}

\begin{abstract}
We formulate several open questions related to enumerative
combinatorics, which arise in the spectral analysis of Hill operators
with trigonometric polynomial potentials.
\end{abstract}

\maketitle


\bigskip

In analysis of asymptotics of spectral gaps and the basisness
property of root function systems of Hill operators
\begin{equation}
\label{1.1}  Ly = - y^{\prime \prime} + v(x)y, \quad   0 \leq  x
\leq \pi,
\end{equation}
with potentials that are trigonometric polynomials we faced (see
\cite{DM11,DM10,DM25a,DM25}) a series of Enumerative Combinatorics
questions. Some of them have been settled, and this allows us to find
the asymptotics of the spectral gap of Schr\"odinger operators with
potentials of the form
\begin{equation}
\label{1.2} v(x)= a \cos 2x + b \cos 4x
\end{equation}
(see \cite[Theorem 7]{DM10}), or to prove (see \cite[Theorem
21]{DM25}) that the Hill operator with potentials of the form
\begin{equation}
\label{1.3} v(x) = a e^{-2ix}+ Ae^{-4ix}+ be^{2ix} +Be^{4ix}, \quad
a,b, A,B \neq 0,
\end{equation}
and periodic (or antiperiodic) boundary conditions has a basis in
$L^2([0,\pi])$ consisting of its root functions if and only if
$|A|=|B|$ (under the restriction that neither $-b^2/(4B) $ nor
$-a^2/(4A)$ is an integer square).

However, many of those questions remain open. They seem to be
interesting by themselves but their answers would imply series of
results on spectra and spectral decompositions related to Hill
operators with trigonometric potentials. Let us formulate some of
them. \bigskip

1. For each $n\in \mathbb{N}$  a {\em walk} $x$ from $-n$ to $n$ or
from $n$ to $-n$ is defined by its {\em sequence of steps}
\begin{equation}
  \label{2.21}
x=(x(t))_{t=1}^{\nu+1}, \quad 1\leq \nu=\nu(x)<\infty,
\end{equation}
where $x(t) \in 2 \mathbb{Z} \setminus \{0\},$ and respectively,
\begin{equation}
  \label{2.22}
\sum_{t=1}^{\nu+1} x(t) = 2n   \quad \text{or}  \quad
\sum_{t=1}^{\nu+1} x(t) = -2n.
\end{equation}
A walk $x$ is called {\em admissible} if its {\em  vertices} $j(t) =
j(t,x)$ given, respectively,  by
\begin{equation}
  \label{2.23}
j(0) = -n  \quad \text{or} \;\;j(0) = +n
\end{equation}
and
\begin{equation}
  \label{2.24}
j(t) =-n + \sum_{i=1}^t x(i) \quad   \text{or}\quad j(t) =
 n + \sum_{i=1}^t x(i), \quad 1\leq t \leq \nu+1,
\end{equation}
satisfy the condition
\begin{equation}
  \label{2.25}
 j(t) \neq \pm n \quad \text{for} \;\; 1\leq t \leq \nu.
\end{equation}
Let $X_n$ and $Y_n$ be, respectively, the set of all admissible
walks from $-n$ to $n$ and from $n$ to $-n, $ and let $\{V(m), \, m
\in 2\mathbb{Z}\}$ be a sequence of complex numbers such that
$V(0)=0. $ For each walk $x\in X_n$ or $x\in Y_n $ we set
\begin{equation}
  \label{2.26}
h_1(x) = \prod_{t=1}^\nu [n^2 - j(t)^2 ]^{-1}, \quad h(x)= h_1(x)
\prod_{t=1}^{\nu+1} V(x(t))
\end{equation}

2. In the sequel we consider finite sets $F$ of permitted steps
$x(t).$ We say that a walk $x$ is $F$-admissible if it is admissible
and $x(t) \in F$ for every $t=1, \ldots, \nu +1. $ We define $X_n^+
(F)$ (or $Y^-_n (F)$) to be the set of all $F$-admissible walks with
positive steps $x(t)>0$ (or negative steps $x(t)<0$).

In these notations, the following statement holds (it is proven in
\cite{DM10}, see formulas (109)-(113) there).

\begin{Proposition}
\label{prop1} Let $F=F_2=\{-2, -4, 2, 4\},$ and let
\begin{equation}
  \label{1.10}
V(-2)=a, \quad V(2) =b, \quad V(-4) =A, \quad V(4) = B.
\end{equation}
Then, for even $n=2m,$
\begin{equation}
  \label{1.11}
\sum_{\xi \in X_n^+} h(\xi) =\frac{1}{4^{m-1}} \prod_{j=1}^m [b^2/4
+ (2j-1)^2 B],
\end{equation}
and for odd $n=2m+1$
\begin{equation}
  \label{1.12}
\sum_{\xi \in X_n^+} h(\xi) =-\frac{b}{4^{m}} \prod_{j=1}^m [b^2/4 +
(2j)^2 B].
\end{equation}
\end{Proposition}

\begin{Question}
\label{q1} Let $F_m =\{j \in 2\mathbb{Z}\setminus \{0\}: \; |j|\leq
m\}, $ and let $z= \{V(2k)\}_1^m \in \mathbb{C}^m. $ Find the
polynomials
\begin{equation}
  \label{1.13}
P_{m,n} (z) = \sum_{x \in X_n^+ (F_M)} h(x)
\end{equation}
and/or their asymptotics for fixed $m$ and $n \to \infty. $
\end{Question}
Proposition \ref{prop1} answers this question for $m=2.$
\bigskip

3. The identities (\ref{1.11}) and (\ref{1.12}) were discovered and
proven by analyzing the combinatorial meaning of the asymptotic
formulas for the spectral gaps of Hill operators with potentials of
the form (\ref{1.2}). From (\ref{1.11}) and (\ref{1.12}) one derives
the following (Theorem 8 in \cite{DM10}).
\begin{Proposition}
\label{prop2}
(a)  If $k$ and  $m, \; 1 \leq k \leq m, $ are fixed,
then
\begin{equation}
\label{4.1} \sum  \prod_{s=1}^k (m^2 -i_s^2) = \sum_{1\leq j_1 <
\cdots < j_k \leq m} \prod_{s=1}^k (2j_s -1)^2,
\end{equation}
where the left sum is over all indices $ i_1, \ldots, i_k $ such
that
$$ -m <i_1 <\cdots <i_k <m ,
\quad |i_s - i_r | \geq 2 \;\;\text{if} \;\; s \neq r.  $$

(b)  If $k$ and  $m, \; 1 \leq k \leq m-1, $ are fixed, then
\begin{equation}
\label{4.2} \sum \prod_{s=1}^k \left [(2m-1)^2 -(2i_s -1)^2 \right]
= \sum_{1\leq j_1 < \cdots < j_k \leq m-1} \prod_{s=1}^k (4j_s)^2,
\end{equation}
where the left sum is over all indices $ i_1, \ldots, i_k $ such
that
$$ -m +1 <i_1 <\cdots <i_k < m ,
\quad |i_s - i_r | \geq 2 \;\;\text{if} \;\; s \neq r.  $$
\end{Proposition}

4. A few years later D. Zagier gave a direct combinatorial
proof\footnote{We have asked many colleagues working in combinatorics
for a direct proof of (\ref{4.1}) and (\ref{4.2}). Finally, such a
proof was given by D. Zagier in August 2010 when we met him by chance
at Mathematisches Forschungsinstitut Oberwolfach.} of (\ref{4.1}) and
(\ref{4.2}). These identities are used in an essential way by J. M.
Borwein, A. Straub, J. Wan and W. Zudilin in \cite{BSWZ11},
\cite{BSWZ12}. Moreover, an appendix to \cite{BSWZ12} presents the
proof of D. Zagier.

Just recently an elegant "elementary" proof was given by
S.~Rosenberg \cite{Ro12}; it is based on information about the
spectrum $\sigma (M_n) = \{n-2\ell: \; 0\leq \ell \leq n\} $ of
Sylvester-Kac matrices (see \cite{TT91} or \cite{EK93}).
\bigskip

5. Good information about the polynomials (\ref{1.13}) would help to
find asymptotics of the sequences
\begin{equation}
\label{9.1} \beta_n^+ = \sum_{x\in X_n (F_m)} h(x), \quad \beta_n^- =
\sum_{y\in Y_n (F_m)} h(y)
\end{equation}
as well because the following holds.
\begin{Proposition}
\label{prop3}
\begin{equation}
\label{9.3} \beta_n^+ = \left ( \sum_{x\in X_n^+} h(x) \right ) \left
[ 1+ O(\log n/n)\right ]
\end{equation}
\begin{equation}
\label{9.4} \beta_n^- = \left ( \sum_{x\in Y_n^+} h(x) \right ) \left
[ 1+ O(\log n/n)\right ]
\end{equation}
if, respectively,
\begin{equation}
\label{9.5}  \sum_{x\in X_n^+} |h(x)| \leq C(V) \left |  \sum_{x\in
X_n^+} h(x)   \right |
\end{equation}
or
\begin{equation}
\label{9.6}  \sum_{x\in Y_n^-} |h(x)| \leq C(V) \left |  \sum_{x\in
Y_n^-} h(x)   \right |
\end{equation}
\end{Proposition}

\begin{Question}
\label{q2} Is it true that (\ref{9.5}) holds for almost all $\{V\}=
z\in \mathbb{C}^m.$  The structure of the exceptional set $W=\{z: \;
C(z) =\infty \}$ is of special interest also.
\end{Question}

6. Suppose $F$ has just two elements, say $F= \{-2R, +2S\}$ and
$V(-2R)=a,$ $V(2S)= b. $  We are concerned on the polynomials
\begin{equation}
\label{10.1} P_\kappa (a,b) =  \sum_{x \in X_n (\kappa)} h(x),
\end{equation}
where $X_n (\kappa)$ is the set of all admissible walks from $-n$ to
$n$ with $\kappa $ negative steps, i.e., $x\in X_n (\kappa) $ if and
only if
\begin{equation}
\label{11.1} \#\{t: \; x(t)= -2R\} =\kappa.
\end{equation}
If
\begin{equation}
\label{11.2} q=\#\{t: \; x(t)= 2S\},
\end{equation}
then $-2Rp + 2Sq = 2n,$ or
\begin{equation}
\label{11.3} Sq = n+R \kappa.
\end{equation}
Of course, if $q= \frac{1}{S}(n+R\kappa) \not \in \mathbb{N},$ then
$X_n (\kappa)= \emptyset $ and $P_\kappa (a,b) \equiv 0. $
Otherwise,
$$
P_\kappa (a,b) =a^\kappa b^q B_\kappa (n), \quad B_\kappa (n) =
\sum_{x \in X_n (\kappa)} h_1 (x)
$$
and we need to analyze just the sequence $B_\kappa (n).$

Different pairs $\{R,S\}$ bring completely different problems
related to solvability of (\ref{11.3}) in $\mathbb{N}.$ We will
mention just one question from this series.

Let $F= \{-2, +4\},$  i.e., $R=1, S=2, $  and let $n= 2m+1.$ Then
(\ref{11.3}) becomes
\begin{equation}
\label{12.2} 2q = 2m+1+ \kappa,
\end{equation}
so $X_n^+ (\kappa) = \emptyset $ if $\kappa $ is even.
Unfortunately,
\begin{equation}
\label{12.3} H_0 = B_1 (2m+1) = \sum_{x\in X_n^+ (1)} h_1 (x) =0
\end{equation}
as we realized in \cite[Proposition 25]{DM25}, see formulas (145) -
(154). (As a matter of fact, (\ref{12.3}) is equivalent to the
fundamental identity for Catalan numbers -- see in \cite{LW92} p.117,
formulas (14.10) - (14.12).)

But the next level $X^+_n (3) $ gives some hopes.

\begin{Question}
\label{q3} Find an asymptotic formula for the sequence $B_3 (2m+1),$
or at least explain that for some absolute constant
\begin{equation}
\label{13.1} |B_3 (2m+1)| \geq C^m/m!.
\end{equation}
\end{Question}
(It is not difficult to see that the upper estimate
\begin{equation}
\label{13.2} |B_3 (2m+1)| \leq \sum_{x \in X_n^+ (3)} |h_1 (x)| \leq
C_1^m/m!
\end{equation}
holds.)

The inequality (\ref{13.1}) would imply that {\em the root function
system of the Hill operator (\ref{1.1}) with potentials of the form
\begin{equation}
\label{14.1} v(x) = a e^{-2ix} + b e^{4ix}, \quad a, b \neq 0,
\end{equation}
subject to antiperiodic boundary conditions does not contain a basis
in} $L^2 ([0,\pi]). $

\end{document}